\documentclass[reqno,a4paper]{amsart}

\usepackage{amsmath,amsfonts,amssymb}
\usepackage{verbatim}
\usepackage{enumerate}
\usepackage{url}
\usepackage{tikz}
\usepackage{setspace}
\usepackage[left=2.5cm,right=2.5cm, top=3cm, bottom=2.5cm]{geometry}
\usepackage[utf8]{inputenc} 
\usepackage[T1]{fontenc}
\usepackage{hyperref}
\usetikzlibrary{arrows}

\theoremstyle{plain}
\newtheorem{theorem}{Theorem}[section]
\newtheorem{lemma}[theorem]{Lemma}

\def\proof{\noindent {\it Proof: }}
\def\qed{\hfill\hbox{$\square$}}

\theoremstyle{definition}

\numberwithin{equation}{section}

\title{Extremal product-one free sequences over $C_n \rtimes_s C_2$}

\author[F. E. Brochero Mart\'{\i}nez]{F. E. Brochero Mart\'{\i}nez}
\address{
Departamento de Matem\'{a}tica\\
Universidade Federal de Minas Gerais\\
UFMG\\
Belo Horizonte, MG\\
31270-901\\
Brazil\\
}
\email{fbrocher@mat.ufmg.br }

\author[S. Ribas]{S. Ribas}
\address{
Departamento de Matem\'{a}tica\\
Universidade Federal de Ouro Preto\\
UFOP\\
Ouro Preto, MG\\
35400-000\\
Brazil\\
}
\email{savio.ribas@ufop.edu.br }

\thanks{The first author was partially supported by FAPEMIG APQ-02973-17, Brazil.}

\date{\today}

\onehalfspace

\begin{document}

\maketitle

\begin{abstract}
Let $G$ be a finite group multiplicatively written. The small Davenport constant of $G$ is the maximum positive integer ${\sf d}(G)$ such that there exists a sequence $S$ of length ${\sf d}(G)$ for which every subsequence of $S$ is product-one free. Let $s^2 \equiv 1 \pmod n$, where $s \not\equiv \pm1 \pmod n$. It has been proven that ${\sf d}(C_n \rtimes_s C_2) = n$ (see \cite[Lemma~6]{ZhGa}). In this paper, we determine all sequences over $C_n \rtimes_s C_2$ of length $n$ which are product-one free. It completes the classification of all product-one free sequences over every group of the form $C_n \rtimes_s C_2$, including the quasidihedral groups and the modular maximal-cyclic groups.
\end{abstract}

\section{Introduction}

Given a finite group $G$ multiplicatively written, the {\em zero-sum problems} study conditions to ensure that a given sequence over $G$ has a non-empty subsequence with prescribed properties (such as  length, repetitions, weights) such that the product of its terms in some order equals $1$.

\subsection{Definitions and notations}

A sequence over a finite group $G$ is a finite and unordered element $S$ of the free abelian monoid $\mathcal F(G)$, equipped with the sequence concatenation product denoted by $\boldsymbol{\cdot}$. The form of $S$ is $$S = g_1 \boldsymbol{\cdot} {\dots} \boldsymbol{\cdot} g_k = \prod_{1 \le i \le k}^{\bullet} g_i = \prod_{1 \le i \le k}^{\bullet} g_{\tau(i)} = \prod_{g \in G}^{\bullet} g^{[v_g(S)]} \in \mathcal F(G),$$ where $g_1, {\dots}, g_k \in G$ are the {\em terms} of $S$, $k = |S| \ge 0$ is the {\em length} of $S$, $\tau: \{1,2,{\dots},k\} \to \{1,2,{\dots},k\}$ is a permutation, and $v_g(S) = \#\{i \in \{1,2,{\dots},k\} \, ; \; g_i = g\}$ is the {\em multiplicity} of $g$ in $S$. We say that a sequence $T$ is a {\em subsequence} of $S$ if $v_g(T) \le v_g(S)$ for all $g \in G$, that is, if $T \mid S$ in $\mathcal F(G)$; in this case, we write $S \boldsymbol{\cdot} T^{[-1]} = \displaystyle\prod_{g \in G}^{\bullet} g^{[v_g(S) - v_g(T)]}$.
We also define:
\begin{align*}
\pi(S) &= \{ g_{\tau(1)} \cdot {\dots} \cdot g_{\tau(k)} \in G; \; \tau \text{ is a permutation of }\{1,2,\dots,k\}\} \quad \text{(the {\em set of products} of $S$)}; \\
\Pi(S) &= \bigcup_{T \mid S \atop |T| \ge 1} \{\pi(T)\} \subset G \quad \text{(the {\em set of subproducts} of $S$)}; \\
S \cap K &= \displaystyle\prod_{{g \mid S \atop g \in K}}^{\bullet} g \quad \text{(the subsequence of $S$ that lie in a subset $K$ of $G$)}. 
\end{align*}
The sequence $S$ is called {\em product-one free} if $1 \not\in \Pi(S)$, and {\em product-one sequence} if $1 \in \pi(S)$.

\subsection{The small Davenport constant}

An important type of zero-sum problem is the {\em small Davenport constant} of $G$. Denoted by ${\sf d}(G)$, this constant is the maximal integer such that there exists a sequence over $G$ (repetition allowed) of length ${\sf d}(G)$ which is product-one free, i.e.,
$${\sf d}(G) = \sup\{|S|>0; \; S \in \mathcal F(G) \text{ is product-one free}\}.$$ 

Denote by $C_n$ the cyclic group of order $n$. It is known that ${\sf d}(C_n) = n-1$. Zhuang \& Gao \cite{ZhGa} showed that ${\sf d}(G) = n$ for every non-abelian group $G$ of the form $C_n \rtimes_s C_2$. The conditions $s^2 \equiv 1 \pmod n$ and $s \not\equiv 1 \pmod n$ are required in order to obtain non-abelian groups. The case $s \equiv -1 \pmod n$ yields $D_{2n}$, the Dihedral Group of order $2n$.

\subsection{Inverse zero-sum problems}

Let $G$ be a finite group. By definition, there exist sequences $S \in \mathcal F(G)$ with ${\sf d}(G)$ terms that are product-one free. The {\em inverse zero-sum problems} study the structure of these extremal length sequences; see \cite{GaGeSc, Sc, GaGeGr} for an overview on inverse problems over abelian groups. In \cite{BEN} (see also Lemma \ref{lemaciclico} below) and \cite[Theorem~3.2]{Sc}, one can find the solution of the inverse problem associated with small Davenport constant over cyclic groups and over abelian groups of rank two, respectively. For non-abelian groups, the first results of this kind have emerged in recent years; see \cite{MR1,MR2,OZ1,OZ2}. In particular, the authors solved in \cite{MR2} the inverse problem associated with small Davenport constant over Dihedral Groups.

Let $n$ be a positive integer such that there exists an integer $s$ satisfying
\begin{equation}\label{hipotese}
s^2 \equiv 1 \!\!\!\pmod n, \quad \text{ but } \quad s \not\equiv \pm1 \!\!\!\pmod n,
\end{equation}
and consider the group $C_n \rtimes_s C_2$ with presentation
\begin{equation}\label{defCnC2}
\big\langle \, x,y \, \mid \, x^2 = y^n = 1, \; yx = xy^s \, \big\rangle, \quad \text{ where the pair $(n,s)$ satisfies the conditions \eqref{hipotese}}.
\end{equation}

The conditions \eqref{hipotese} are necessary, since the semidirect product $C_n \rtimes_s C_2$ should yield a non-abelian group. Indeed, $s \equiv 1 \pmod n$ means that $C_n \rtimes_s C_2$ is an abelian group ($C_{2n}$ or $C_2 \times C_n$, depending on whether $n$ is odd or even, respectively), while $s \equiv -1 \pmod n$ means that $C_n \rtimes_s C_2$ is the Dihedral Group of order $2n$, for which the inverse problems are already solved. It is worth mentioning that this family of metacyclic groups includes the quasidihedral groups, when $(n,s) = (2^t,2^{t-1}-1)$, $t \ge 3$, and the modular maximal-cyclic groups, when $(n,s) = (2^t,2^{t-1}+1)$, $t \ge 3$.

In this paper, we solve the inverse problem over $C_n \rtimes_s C_2$ related to small Davenport constant, proving the following result:

\begin{theorem}\label{thminverse}
Let $G \simeq C_n \rtimes_s C_2$ be the group with presentation given by \eqref{defCnC2}. Let $S \in \mathcal F(G)$ be a sequence of length $|S| = n$. Then $S$ is product-one free if and only if 
\begin{enumerate}[(i)]
\item For any $(n,s)$, $S = (y^u)^{[n-1]} \boldsymbol{\cdot} xy^v \quad \text{ for some integers $u,v$ with } \gcd(u,n) = 1$, or
\item For $(n,s) = (2^t,2^{t-1}+1)$ where $t \ge 3$, 
$S = 
\begin{cases}
(xy^u)^{[n-1]} \boldsymbol{\cdot} xy^v &\text{for $u$ odd and $v$ even, or} \\
(xy^u)^{[n-1]} \boldsymbol{\cdot} y^v &\text{for $u, v$ odd}.
\end{cases}$
\end{enumerate}
\end{theorem}

In particular, together with \cite[Theorem~1.3]{MR2} and \cite[Theorem~3.2]{Sc}, the previous theorem completes the classification of all the extremal length product-one free sequences over all groups of the form $C_n \rtimes_s C_2$. The case {\it (ii)} of previous theorem appears to have some extra extremal product-one free sequences, but this is not actually true. It turns out that $C_n \rtimes_s C_2 = \langle x,xy^v \mid v \text{ is odd} \rangle$ for $(n,s) = (2^t, 2^{t-1}+1)$ and $t \ge 3$, that is, the extra sequences arise only because of this change of generators.

The paper is organized as follows: In Section \ref{lemas}, we present two lemmas which will be used later. 
In Section \ref{proofmain}, we prove the main result by splitting the proof into cases according to the number of terms of $S$ that belong to the normal subgroup of order $n$, and in some cases we further split into subcases according to either the value of $s$ or the factorization of $n$.

\section{Preliminary results}\label{lemas}

In this section, we present two preliminary results that will be used in the proof of the main theorem. The first can be seem as a generalization of the inverse zero-sum problem over cyclic groups, dealing with product-one free sequences of large length and their subsequence sums.

\begin{lemma}[See {\cite[Theorem~11.1]{Gr1}}]\label{lemaciclico}
Let $m \ge 3$ be an integer and fix a generator $y$ of $C_m$. Let $S = \prod_{1 \le i \le |S|}^{\bullet} y^{a_i} \in \mathcal F(C_m)$ be a product-one free sequence of length $|S| > m/2$. Then there exists $y^{a_i} \mid S$ such that $v_{y^{a_i}}(S) \ge \max\left\{ m - 2|S| + 1 , |S| - \left\lfloor \frac{m-1}{3} \right\rfloor \right\}$. In addition, for a residue class $r \pmod m$, denote by $\overline{r}$ the integer such that $0 \le \overline{r} \le m-1$ and $r \equiv \overline{r} \pmod m$. Then there exists an integer $t$ with $\gcd(t,m) = 1$ and $\sum_{1 \le i \le |S|} \overline{a_it} < m$. Moreover, for every $1 \le k \le \sum_{1 \le i \le |S|} \overline{a_it}$, there exists $T_k \mid S$ such that $\sum_{y^{a_i} \mid T_k} \overline{a_it} = k$.
\end{lemma}

The second lemma is related to the factorization of $n$ and has essentially been considered in \cite[Section 5]{GeGr}, nevertheless we summarize it here for convenience. We notice that conditions \eqref{hipotese} guarantee that $n$ can be neither an odd prime power nor twice an odd prime power, otherwise we would have $s \equiv \pm1 \pmod n$. Assuming that $n$ is not a power of $2$, we are going to show that these conditions suffice to factor $n$ in a helpful way.

\begin{lemma}\label{lemma}
Let $n \ge 8$ and $s$ be positive integers satisfying the conditions \eqref{hipotese}.
\begin{itemize}
\item If both $n \neq p^t$ and $n \neq 2p^t$ for every prime $p$ and every integer $t \ge 1$, then there exist coprime integers $n_1, n_2 \ge 3$ such that $s \equiv -1 \pmod {n_1}$, $s \equiv 1 \pmod {n_2}$, and either (A) $n = n_1n_2$ or (B) $n = 2n_1n_2$.
\item If $n = 2^t$ for some $t \ge 3$, then (B) $n = 2n_1n_2$, where either $(n_1,n_2) = (1,2^{t-1})$ satisfies $s \equiv 1 \pmod {n_2}$ or $(n_1,n_2) = (2^{t-1},1)$ satisfies $s \equiv -1 \pmod {n_1}$.
\end{itemize}
\end{lemma}

\proof
Let $n = 2^tm$, where $m$ is odd and $t \ge 0$ is an integer. Since $m$ divides $s^2-1$ and $\gcd(s-1,s+1) \in \{1,2\}$ (depending on $s$ is even or odd), each prime power factor of $m$ divides either $s+1$ or $s-1$. Let $m_1 = \gcd(m,s+1)$ and $m_2 = \gcd(m,s-1)$, so that $m = m_1m_2$. In addition, $s^2 \equiv 1 \pmod {2^t}$ implies that either $s \equiv \pm1 \pmod {2^t}$ or $s \equiv 2^{t-1} \pm 1 \pmod {2^t}$ (whether $t \ge 3$). We consider some cases:
\begin{enumerate}[(i)]
\item {\bf CASE $t = 0$.} In this case, $n = m = m_1m_2$. We set $n_1 = m_1$ and $n_2 = m_2$, hence $n = n_1n_2$ is the desired factorization, as in (A). This is the only case where $s$ can be even; in the following, $s$ must be odd.
\item {\bf CASE $t = 1$.} It is possible to set either $n_1 = 2m_1$ and $n_2 = m_2$ or $n_1 = m_1$ and $n_2 = 2m_2$, hence $n = n_1n_2$ is the desired factorization, as in (A).
\item {\bf CASE $t \ge 2$ and $m \ge 3$.} If $s \equiv -1 \pmod {2^t}$, then we set $n_1 = 2^tm_1$ and $n_2 = m_2$. If $s \equiv 1 \pmod {2^t}$, then we set $n_1 = m_1$ and $n_2 = 2^tm_2$. Therefore, $n = n_1n_2$ is the desired factorization, as in (A). \\
In the case that $t \ge 3$, it is possible that $s \equiv 2^{t-1} \pm 1 \pmod {2^t}$. If $s \equiv 2^{t-1}-1 \pmod {2^t}$, then we set $n_1 = 2^{t-1}m_1$ and $n_2 = m_2$. If $s \equiv 2^{t-1}+1 \pmod {2^t}$, then we set $n_1 = m_1$ and $n_2 = 2^{t-1}m_2$. Thus, $n = 2n_1n_2$ is the desired factorization as in (B).
\item {\bf CASE $t \ge 3$ and $m=1$.} In this case, $s \equiv 2^{t-1} \pm 1 \pmod {2^t}$, which implies that $s \equiv \pm1 \pmod {2^{t-1}}$. For the negative sign, it follows that $(n_1,n_2) = (2^{t-1},1)$, and for the positive sign, it follows that $(n_1,n_2) = (1,2^{t-1})$. Therefore, $n = 2 n_1 n_2$ is the factorization as in (B).
\end{enumerate}
\qed

We highlight that the case (iv) of previous lemma yields the quasidihedral group provided $s \equiv 2^{t-1}-1 \pmod {2^t}$ and the modular maximal-cyclic group provided $s \equiv 2^{t-1}+1 \pmod {2^t}$.

From the previous lemma and Chinese Remainder Theorem, there exists a natural projection 
\begin{equation}\label{isomorphism}
\Psi: C_n \to C_{n_1} \oplus C_{n_2}
\end{equation}
satisfying $\Psi(e) = (e_1,e_2)$ and $\Psi(s \cdot e) = (-e_1,e_2)$, where $C_{n_1} = \langle e_1 \rangle$, $C_{n_2} = \langle e_2 \rangle$ and $C_n = \langle e \rangle$. In case (A) of previous lemma, $\Psi$ is an isomorphism. 

\section{Proof of Theorem \ref{thminverse}}\label{proofmain}

The sequences of the form $(y^u)^{[n-1]} \boldsymbol{\cdot} xy^v \in \mathcal F(C_n \rtimes_s C_2)$, where $\gcd(u,n) = 1$, are product-one free. In addition, if $(n,s) = (2^t,2^{t-1}+1)$ for some $t \ge 3$, then the sequences $(xy^u)^{[n-1]} \boldsymbol{\cdot} xy^v \in \mathcal F(C_n \rtimes_s C_2)$, for $u$ odd and $v$ even, and $(xy^u)^{[n-1]} \boldsymbol{\cdot} y^v \in \mathcal F(C_n \rtimes_s C_2)$, for $u, v$ odd, are product-one free. In fact: 
\begin{itemize}
\item 
If $(xy^u)^a = 1$ for some $a \in \{1,2,\dots,n-1\}$, then $a$ must be even, thus 
$$1 = (y^{2u+n/2})^{a/2} = y^{(2u+n/2)a/2} = y^{(u + n/4)a},$$ 
which implies that $n = 2^t$ divides $a(u+n/4)$. Since $\gcd(u+n/4,n) = 1$, $n$ divides $a$, a contradiction.

\item 
If $(xy^u)^a (xy^v) = 1$, then $a$ must be odd, thus 
$$\quad\quad 1 = [(xy^u)^2]^{(a-1)/2} \cdot (xy^u \cdot xy^v) = y^{(2u+n/2)(a-1)/2} \cdot y^{u+n/2+v} = y^{(2u+n/2)(a-1)/2 + u+v + n/2}.$$
It implies that $n=2^t$ divides $(2u+n/2)(a-1)/2 + u+v + n/2$, but this is a contradiction since $(2u + n/2)(a-1)/2 + u + v + n/2$ is odd.

\item 
If $(xy^u)^a y^v = 1$, then $a$ must be even, thus 
$$\quad 1 = [(xy^u)^2]^{a/2} \cdot y^u = y^{(2u+n/2)a/2} \cdot y^v = y^{(2u+n/2)a + v}.$$
It implies that $n=2^t$ divides $(2u+n/2)a/2 + v$, but this is a contradiction since $(2u + n/2)a/2 + v$ is odd.
\end{itemize}

Therefore it only remains to verify the converse. Let $S \in \mathcal F(C_n \rtimes_s C_2)$ be a product-one free sequence of length $|S| = n$. Let $H = \langle y \rangle \simeq C_n$ be the normal subgroup of $C_n \rtimes_s C_2$ and let $N = xH$. We split into some cases and subcases:
\begin{enumerate}[(a)]
\item {\bf CASE $\mathbf{|S \cap H| = n}$.} In this case, $S$ is contained in the cyclic subgroup of order $n$. Since ${\sf d}(C_n) = n-1$, $S$ cannot be product-one free.
\vspace{1mm}

\item {\bf CASE $\mathbf{|S \cap H| = n-1}$.} In this case, by Lemma \ref{lemaciclico}, the terms of $S \cap H$ must all be equal to the same fixed generator of $H$. Therefore, there exist integers $t$ and $u$ with $\gcd(t,n)=1$ such that $S = (y^t)^{[n-1]} \boldsymbol{\cdot} xy^u$.
\vspace{1mm}

\item {\bf CASE $\mathbf{|S \cap H| = n-k}$, where $\mathbf{2 \le k < n/2}$.} By Lemma \ref{lemaciclico}, we may assume without loss of generality that $S \cap H = \prod_{1 \le i \le n-k}^{\bullet} y^{a_i}$, where $n-k \le \sum_{1 \le i \le n-k} a_i \le n-1$. Moreover, Lemma \ref{lemaciclico} also ensures that $y^{\gamma} \in \Pi(S \cap H)$ for every $1 \le \gamma \le n-k$. If there exists a subsequence $xy^{\alpha} \boldsymbol{\cdot} xy^{\beta} \mid S \cap N$ such that $\alpha s + \beta$ has a representant in the range $[k, n-1]$ modulo $n$, then there exists $y^{\gamma} \in \Pi(S \cap H)$, $1 \le \gamma \le n-k$, such that $xy^{\alpha} \cdot xy^{\beta} \cdot y^{\gamma} = 1$. Therefore, we assume that $\alpha s + \beta$ has a representant in $[1, k-1]$ for all $\alpha, \beta$ such that $xy^{\alpha} \boldsymbol{\cdot} xy^{\beta} \mid S$. Since $s^2 \equiv 1 \pmod n$, we may assume that $s$ has an integer representant in the range $[\sqrt{n+1}, n-\sqrt{n+1}]$. The elements from $S \cap N$ are paired forming more $\lfloor k/2 \rfloor$ terms in $H$, say $y^{a_i}$ for $n-k+1 \le i \le n-k+\lfloor k/2 \rfloor$, and it yields a new sequence $T \in \mathcal F(H)$ of length $|T| = n-k + \lfloor k/2 \rfloor \ge 3n/4$, that is, $T = \prod_{1 \le i \le n-k+\lfloor k/2 \rfloor}^{\bullet} y^{a_i}$ is the concatenation of $S \cap H$ and the pairs of $S \cap N$. If $T$ is product-one free, then the average of the exponents of the terms other than $y$ is at least two, therefore 
$$\frac{n-1 - v_y(T)}{3n/4 - v_y(T)} \ge \frac{\sum_{1 \le i \le |T|} a_i - v_y(T)}{n - k + \lfloor k/2 \rfloor - v_y(T)} \ge 2,$$
which implies $v_y(T) \ge \lceil n/2 \rceil + 1$. 

We are going to show that at least one of the $y \mid T$ comes from a pair of $S \cap N$, that is, $v_y(T) > v_y(S \cap H) = v_y(S)$. Otherwise, the average of the $\lfloor k/2 \rfloor$ terms of $T$ coming from $S \cap N$ is at most $$\frac{\sum_{y^{a_i} \mid T} a_i - \sum_{y^{a_i} \mid S \cap H} a_i}{\lfloor k/2 \rfloor} \le \frac{n-1-(n-k)}{\lfloor k/2 \rfloor} = 
\begin{cases}
\frac{2(k-1)}{k} < 2 \quad &\text{ if $k$ is even,} \\
2 \quad &\text{ if $k$ is odd.}
\end{cases}$$
The Pigeonhole Principle implies that $\Pi(S \cap N) = \{y^2\}$ and $k$ is odd. Let $xy^{\alpha} \boldsymbol{\cdot} xy^{\beta} \mid S \cap N$ such that $xy^{\alpha} \cdot xy^{\beta} = y^2$, thus $xy^{\beta} \cdot xy^{\alpha} = y^{2s}$. In general, it increases the exponent and contradicts the average at most two. In fact, if $2s \equiv 1 \pmod n$, then we would be done since $y^{2s} = y \mid T$. Moreover, $2s \equiv 2 \pmod n$ if and only if either $s \equiv 1 \pmod n$ if $n$ is odd, a contradiction, or $s \equiv n/2+1 \pmod n$ if $n$ is even. In the case that $n$ is even, it follows that $n$ is divisible by $4$ since $\gcd(s,n)=1$. In addition, $S \cap N = (xy^{\alpha})^{[k]}$, where $\alpha s + \alpha \equiv 2 \pmod n$. In fact, if $k \ge 3$, then let $xy^{\alpha} \boldsymbol{\cdot} xy^{\beta} \boldsymbol{\cdot} xy^{\gamma} \mid S \cap N$. It follows that $\alpha s + \beta \equiv 2 \pmod n$ and $\alpha s + \gamma \equiv 2 \pmod n$, from where we obtain $\beta \equiv \gamma \pmod n$ and similarly we obtain $S \cap N = (xy^{\alpha})^{[k]}$. If $k=2$, then $S \in \{y^{[n-2]} \boldsymbol{\cdot} xy^{\alpha} \boldsymbol{\cdot} xy^{\beta} \; , \; y^{[n-3]} \boldsymbol{\cdot} y^2 \boldsymbol{\cdot} xy^{\alpha} \boldsymbol{\cdot} xy^{\beta} \}$, and in both cases it is easy to verify that $S$ is not product-one free since either $\alpha s + \beta \not\equiv 1 \pmod n$ or $\beta s + \alpha \not\equiv 1 \pmod n$. Therefore, $S \cap N = (xy^{\alpha})^{[k]}$. If $\alpha$ is even, then $\alpha s + \alpha \equiv 2 \pmod n$ and $s \equiv n/2 + 1 \pmod n$ imply that $2\alpha \equiv 2 \pmod n$, therefore $\alpha \equiv 1 \pmod {n/2}$ and $\alpha$ is odd, contradiction. If $\alpha$ is odd, then $n/2 + 2\alpha \equiv 2 \pmod n$, therefore $\alpha \equiv n/4 + 1 \pmod {n/2}$ and hence $\alpha \in \{n/4 + 1, 3n/4 + 1\}$. For each value of $\alpha$, we obtain a contradiction using the following products:
$$\begin{cases}
xy^{n/4+1} \cdot y \cdot xy^{n/4+1} \cdot y^{n/2 - 3} = 1 \quad &\text{ if $\alpha = n/4+1$ and $n \equiv 0 \pmod 8$}, \\
xy^{n/4+1} \cdot y^2 \cdot xy^{n/4+1} \cdot y^{n/2 - 4} = 1 \quad &\text{ if $\alpha = n/4+1$ and $n \equiv 4 \pmod 8$}, \\
xy^{3n/4+1} \cdot y \cdot xy^{3n/4+1} \cdot y^{n/2 - 3} = 1 \quad &\text{ if $\alpha = 3n/4+1$ and $n \equiv 0 \pmod 8$}, \\
xy^{3n/4+1} \cdot y^2 \cdot xy^{3n/4+1} \cdot y^{n/2 - 4} = 1 \quad &\text{ if $\alpha = 3n/4+1$ and $n \equiv 4 \pmod 8$},
\end{cases}$$
which implies that $v_y(T) > v_y(S)$ so far.

Let $xy^{\alpha} \boldsymbol{\cdot} xy^{\beta} \mid S \cap N$ with $y = xy^{\alpha} \cdot xy^{\beta}$, so that $y^s = xy^{\beta} \cdot xy^{\alpha}$. We consider the following subcases:

\begin{enumerate}[{(c.}1{)}]
\item {\bf Subcase $n/2 < s \le n-\sqrt{n+1}]$.} Since $1 \le n-s < n/2$ and $|S \cap H| > n/2$, Lemma \ref{lemaciclico} ensures that $y^{n-s} \in \Pi(S \cap H)$. Therefore $xy^{\beta} \cdot xy^{\alpha} \cdot y^{n-s} = 1$, and we are done.

\item {\bf Subcase $\sqrt{n+1} \le s < n/2$.} We proceed similarly to the baby-step giant-step method for the calculus of discrete logarithm \cite{Sh}. Let $\ell \ge 2$ be the greatest integer such that $\ell s \le n$. It follows that $\ell \le \left\lfloor \frac{n}{\sqrt{n+1}} \right\rfloor = \lfloor \sqrt{n} \rfloor$ and $n - s < \ell s \le n$. Since $v_y(T) \ge \lceil n/2 \rceil + 1$, we consider the following product
\begin{equation}\label{eqls}
xy^{\beta} \cdot \underbrace{y \cdot {\dots} \cdot y}_{\ell - 1 \text{ times}} \cdot xy^{\alpha} = y^{\ell s}.
\end{equation}
If $\sqrt{n+1} \le s \le \lceil n/2 \rceil + 1 - \lfloor \sqrt{n} \rfloor$, then $\ell + s \le \lceil n/2 \rceil + 1$ and we complete the product in \eqref{eqls} as $y^{\ell s} \cdot \underbrace{y \cdot {\dots} \cdot y}_{n - \ell s \text{ times}} = 1$ since $1 \le n - \ell s < s$.

Otherwise, $\lceil n/2 \rceil - \lfloor \sqrt{n} \rfloor + 2 \le s < n/2$ implies that $\ell = 2$, therefore $2 + s \le \lceil n/2 \rceil + 1$ and in the same way we obtain $y^{2s} \cdot y^{n - 2s} = 1$, since $n-2s \le \lceil n/2 \rceil - 1$. Hence we are done.
\end{enumerate}

\vspace{1mm}

\item {\bf CASE $\mathbf{|S \cap H| = n-k}$, where $\mathbf{n/2 \le k \le n}$.} 
\begin{enumerate}[({d.}1)]
\item {\bf Subcase} $n = n_1n_2$. 
Group the terms of $S \cap N = xy^{\beta_1} \boldsymbol{\cdot} {\dots} \boldsymbol{\cdot} xy^{\beta_k}$ in pairs $xy^{\beta_i} \boldsymbol{\cdot} xy^{\beta_j}$ in such way that $\beta_i s + \beta_j \equiv 0 \pmod {n_1}$. This is possible always that $\beta_i \equiv \beta_j \pmod {n_1}$ because $s \equiv -1 \pmod {n_1}$, hence we group these pairs until remains $\ell \le n_1$ terms, each of them having the exponents of $y$ in distinct residue classes modulo $n_1$. The number of pairs formed is $\lceil \frac{k-\ell}{2} \rceil$. The $\ell$ remaining terms are grouped into $\lfloor \ell/2 \rfloor$ pairs with product in $H$. We also group these $\lfloor \ell/2 \rfloor$ terms and the terms of $S \cap H = y^{\alpha_1} \boldsymbol{\cdot} {\dots} \boldsymbol{\cdot} y^{\alpha_{n-k}}$ into $\lfloor \frac{n-k + \lfloor \ell/2 \rfloor}{n_1} \rfloor$ subsequences, each of them is not product-one free in the quotient $C_n / C_{n_1} \simeq \langle y^{n_1} \rangle$. Thus, we obtain $\frac{k-\ell}{2} + \lfloor \frac{n-k + \lfloor \ell/2 \rfloor}{n_1} \rfloor$ elements in $\langle y^{n_1} \rangle \simeq C_{n_2}$. Since $\ell \le n_1$, we have that
\begin{align*}
\frac{k-\ell}{2} + \left\lfloor \frac{n-k + \lfloor \ell/2 \rfloor}{n_1} \right\rfloor &\ge \frac{k-\ell}{2} + \frac{n-k + \frac{\ell-1}{2} - (n_1 - 1)}{n_1} \\
&\ge \frac{k}{2} - \frac{n_1}{2} + n_2 - \frac{k}{n_1} - \frac{1}{2} + \frac{1}{2n_1},
\end{align*}
and if the latter is bigger than ${\sf d}(C_{n_2}) = n_2 - 1$, then we are done, since we obtain more than ${\sf d}(C_{n_2})$ elements in $\langle y^{n_1} \rangle$, which yields a product-one sequence over $C_n$. But this assertion is equivalent to $k > n_1 + 1 + \frac{3}{n_1 - 2}$. Since $k \ge \frac{n}{2} = \frac{n_1n_2}{2}$, if $\frac{n_1n_2}{2} > n_1+1+\frac{3}{n_1-2}$, then we are done. It implies that we are done unless either $(n_1,n_2) \in (3,4)$ and $6 \le k \le 7$ or $(n_1,n_2) = (4,3)$ and $k = 6$, since otherwise $k > n_1+1+\frac{3}{n_1-2}$. But these cases follows directly from the inequality $$\frac{k-\ell}{2} + \left\lfloor \frac{n-k + \lfloor \ell/2 \rfloor}{n_1} \right\rfloor \ge n_2.$$

\item {\bf Subcase} $n = 2n_1n_2$, $n_1$ even and $n_2 \ge 3$.
We slightly modify the previous argument to obtain $\big\lfloor \frac{n-k}{2n_1} \big\rfloor + \big\lfloor \frac{k-2n_1}{2} \big\rfloor$ elements in $\langle y^{2n_1} \rangle \simeq C_n/C_{2n_1} \simeq C_{n_2}$. If
$$\left\lfloor \frac{n-k}{2n_1} \right\rfloor + \left\lfloor \frac{k-2n_1}{2} \right\rfloor \ge \frac{n-k - 2n_1+1}{2n_1} + \frac{k - 2n_1 - 1}{2} > n_2 - 1 = {\sf d}(C_{n_2}),$$ 
then we are done, since in this case we obtain more than ${\sf d}(C_{n_2})$ elements in $\langle y^{2n_1} \rangle$ which yields a product-one subsequence over $C_n$. But the last inequality above is equivalent to $k > 2n_1 + 3+ \frac{2}{n_1-1}$.

Since $k \ge \frac{n}{2} = n_1n_2$, if $n_1n_2 > 2n_1+3+\frac{2}{n_1-1}$, then we are done. Since $n_1 \ge 4$ and $n_2 \ge 3$, it follows that $n_1n_2 \ge 3n_1 > 2n_1 + 4 \ge 2n_1 + 3 + \frac{2}{n_1-1}$, hence we are done.

\item {\bf Subcase} $n = 2n_1n_2$, $n_2$ even and $n_1 \ge 3$.
Similar to the previous case, we obtain $\big\lfloor \frac{n-k}{n_1} \big\rfloor + \big\lfloor \frac{k-n_1}{2} \big\rfloor$ elements in $\langle y^{n_1} \rangle \simeq C_n/C_{n_1} \simeq C_{2n_2}$. If
$$\left\lfloor \frac{n-k}{n_1} \right\rfloor + \left\lfloor \frac{k-n_1}{2} \right\rfloor \ge \frac{n-k - n_1+1}{n_1} + \frac{k - n_1 - 1}{2} > 2n_2 - 1 = {\sf d}(C_{2n_2}),$$ 
then we are done, since in this case we obtain more than ${\sf d}(C_{2n_2})$ elements in $\langle y^{n_1} \rangle$ which yields a product-one subsequence over $C_n$. But the last inequality above is equivalent to $k > n_1 + 3+ \frac{4}{n_1-2}$.

Since $k \ge \frac{n}{2} = n_1n_2$, if $n_1n_2 > n_1+3+\frac{4}{n_1-2}$, then we are done. Since $n_1 \ge 3$ and $n_2 \ge 4$, it follows that $n_1n_2 \ge 4n_1 > n_1 + 7 \ge n_1 + 3 + \frac{4}{n_1-2}$, hence we are done.

\item {\bf Subcase} $n = 2^t = 2n_1n_2$, $n_2 = 1$. 
In this case, $s \equiv n_1 - 1 \pmod n$, thus $s \equiv -1 \pmod {n_1}$. Let $S = \prod_{1 \le i \le k}^{\bullet} (xy^{\alpha_i}) \boldsymbol{\cdot} \prod_{1 \le i \le n-k}^{\bullet} (y^{\beta_i})$. We are going to look at the $\alpha_i$'s modulo $n/2$ and, for this, we consider the following values of $k$:

\begin{enumerate}[{(d.4.}i{)}]
\item Sub-subcase $k \ge n/2 + 3$. Since $k - 2 \ge n/2+1$, there exist exponents (say) $\alpha_1 \equiv \alpha_2 \pmod {n/2}$ and $\alpha_3 \equiv \alpha_4 \pmod {n/2}$. It implies that either $$\quad \quad \quad \quad \quad \quad xy^{\alpha_1} \cdot xy^{\alpha_2} = 1 \quad \text{or} \quad xy^{\alpha_3} \cdot xy^{\alpha_4} = 1 \quad \text{or} \quad xy^{\alpha_1} \cdot xy^{\alpha_2} \cdot xy^{\alpha_3} \cdot xy^{\alpha_4} = 1.$$
In fact, the previous first two products are $y^{\alpha_2 - \alpha_1} \in \{1,y^{n/2}\}$ and $y^{\alpha_4 - \alpha_3} \in \{1,y^{n/2}\}$. If both of them are distinct than $1$, then the third product is $y^{\alpha_2 - \alpha_1 + \alpha_4 - \alpha_3} = 1$, and we are done.

\item Sub-subcase $k = n/2 + 2$. Since $k - 2 = n/2$, there exist exponents (say) $\alpha_1 \equiv \alpha_2 \pmod {n/2}$ and either [$\alpha_3 \equiv \alpha_4 \pmod {n/2}$] or [$\alpha_{i+3} \equiv i \pmod {n/2}$ for $0 \le i \le n/2-1$]. If $\alpha_3 \equiv \alpha_4 \pmod {n/2}$, then we use the same argument than Sub-subcase (d.4.i). Otherwise, $\alpha_{i+3} \equiv i \pmod {n/2}$ for $0 \le i \le n/2-1$ implies that 
$$xy^{\alpha_{i+3}} \cdot xy^{\alpha_{i+5}} \cdot xy^{\alpha_{i+6}} \cdot xy^{\alpha_{i+4}} \in \{1,y^{n/2}\},$$
therefore either 
\begin{align*}
xy^{\alpha_1} \cdot xy^{\alpha_2} &= 1 \quad \text{or} \\
xy^{\alpha_{i+3}} \cdot xy^{\alpha_{i+5}} \cdot xy^{\alpha_{i+6}} \cdot xy^{\alpha_{i+4}} &= 1 \quad \text{or} \\
xy^{\alpha_1} \cdot xy^{\alpha_2} \cdot xy^{\alpha_{i+3}} \cdot xy^{\alpha_{i+5}} \cdot xy^{\alpha_{i+6}} \cdot xy^{\alpha_{i+4}} &= 1,
\end{align*}
and we are done.

\item Sub-subcase $k = n/2 + 1$. Since $k - 2 = n/2 - 1$, there exist exponents (say) $\alpha_1 \equiv \alpha_2 \pmod {n/2}$ and either [$\alpha_3 \equiv \alpha_4 \pmod {n/2}$] or [there exist $0 \le \ell \le n/2-1$ such that the sets $\{\alpha_{i+3}; 0 \le i \le n/2-2\}$ and $\{0,1,2,\dots,n/2-1\}\setminus \{\ell\}$ are the same modulo $n/2$]. If $\alpha_3 \equiv \alpha_4 \pmod {n/2}$, then we use the same argument than Sub-subcase (d.4.i). Otherwise, there exist $0 \le \ell \le n/2-1$ such that the sets $\{\alpha_{i+3}; 0 \le i \le n/2-2\}$ and $\{0,1,2,\dots,n/2-1\}\setminus \{\ell\}$ are the same modulo $n/2$. If $t \ge 4$ (that is, $n \ge 16$), then the set $\{\alpha_{i+3}; 0 \le i \le n/2-2\}$ has four consecutive elements, therefore we can apply the same argument as Subcase (d.4.ii). If $t=3$ (that is, $n = 8$), then $|S \cap H| = n-k = 3$ and $xy^{\alpha_3} \cdot xy^{\alpha_4}$ generates another product in $C_n = C_8$. These four elements in $C_8$ yields a subsequence $T \mid (S \cap H) \boldsymbol{\cdot} xy^{\alpha_3} \boldsymbol{\cdot} xy^{\alpha_4}$ such that either $1 \in \pi(T)$ or $y^4 \in \pi(T)$. Therefore either 
$$\quad \quad \quad xy^{\alpha_1} \cdot xy^{\alpha_2} = 1 \quad \text{or} \quad 1 \in \pi(T) \quad \text{or} \quad 1 \in \pi(T \boldsymbol{\cdot} xy^{\alpha_1} \boldsymbol{\cdot} xy^{\alpha_2}),$$
thus we are done.

\item Sub-subcase $k = n/2$. Since $|S \cap H| = n/2$, there exists $T \mid S \cap H$ such that either $1 \in \pi(T)$ or $y^{n/2} \in \pi(T)$. Since $|S \cap N| = n/2$, either there exist exponents (say) $\alpha_1 \equiv \alpha_2 \pmod {n/2}$ or [$\alpha_{i+1} \equiv i \pmod {n/2}$ for $0 \le i \le n/2-1$]. If $\alpha_1 \equiv \alpha_2 \pmod {n/2}$, then either 
$$\quad \quad \quad xy^{\alpha_1} \cdot xy^{\alpha_2} = 1 \quad \text{or} \quad 1 \in \pi(T) \quad \text{or} \quad 1 \in \pi(T \boldsymbol{\cdot} xy^{\alpha_1} \boldsymbol{\cdot} xy^{\alpha_2}),$$
thus we are done. Otherwise, $\alpha_{i+1} \equiv i \pmod {n/2}$ for $0 \le i \le n/2-1$ implies that $xy^{\alpha_1} \cdot xy^{\alpha_3} \cdot xy^{\alpha_4} \cdot xy^{\alpha_2} \in \{1,y^{n/2}\}$. Therefore either 
$$\quad \quad \quad \quad \quad \quad \quad 1 \in \pi(T) \quad \text{or} \quad xy^{\alpha_1} \cdot xy^{\alpha_3} \cdot xy^{\alpha_4} \cdot xy^{\alpha_2} = 1 \quad \text{or} \quad 1 \in \pi \left( T \boldsymbol{\cdot} \prod_{1 \le i \le 4}^{\bullet} (xy^{\alpha_i}) \right),$$
hence we are done.
\end{enumerate}

\item {\bf Subcase} $n = 2^t = 2n_1n_2$, $n_1 = 1$. 
In this case, $s \equiv n_2 + 1 \pmod n$. Let $$S = S_1 \boldsymbol{\cdot} S_2 \boldsymbol{\cdot} S_3 \boldsymbol{\cdot} S_4,$$ where 
$$
S_1 = \prod_{y^{\alpha} \mid S \atop \alpha \text{ odd}}^{\bullet} y^{\alpha}, \quad
S_2 = \prod_{y^{\alpha} \mid S \atop \alpha \text{ even}}^{\bullet} y^{\alpha}, \quad
S_3 = \prod_{xy^{\alpha} \mid S \atop \alpha \text{ odd}}^{\bullet} xy^{\alpha}, \quad
S_4 = \prod_{xy^{\alpha} \mid S \atop \alpha \text{ even}}^{\bullet} xy^{\alpha}. 
$$

We have that
$$xy^{\alpha} \cdot xy^{\beta} = y^{\alpha (n/2+1) + \beta} = 
\begin{cases}
y^{\alpha + \beta} & \text{ if $\alpha$ is even,} \\
y^{\alpha + \beta + n/2} & \text{ if $\alpha$ is odd.}
\end{cases}$$

Let $C_{n/2} \simeq \langle y^2 \rangle$. We see that $S_2 \in \mathcal F(C_{n/2})$. Furthermore, for $i \in \{1,3,4\}$, the product of any two terms of $S_i$ belongs to $C_{n/2}$. Therefore the terms can be grouped into pairs whose products belong to $C_{n/2}$. In this way, we obtain 
$$L = \left\lfloor \frac{|S_1|}{2} \right\rfloor + |S_2| + \left\lfloor \frac{|S_3|}{2} \right\rfloor + \left\lfloor \frac{|S_4|}{2} \right\rfloor \ge \frac{n + |S_2| - 3}{2} \ge \frac{n-3}{2}$$
disjoint products in $C_{n/2}$. If $L > {\sf d}(C_{n/2}) = n/2 - 1$, then $S$ is not product-one free. Therefore we may assume that $L = n/2 - 1$. In particular, $|S_2| \le 1$. If $|S_2| = 1$ and $L = n/2 - 1$, then $|S_1|$, $|S_3|$ and $|S_4|$ are all odd. It implies that there is one term left for each $S_1$, $S_3$ and $S_4$, and their product (in any order) lies in $C_{n/2}$, hence we obtain one element more in $C_{n/2}$.

From now on, we assume that $|S_2| = 0$. By Lemma \ref{lemaciclico}, these $n/2 - 1$ terms of $C_{n/2}$ must comprise the same fixed generator of $C_{n/2}$, say equal to $y^2$ without loss of generality.

Suppose that $|S_3| \ge 3$ and let $xy^{\alpha_1} \boldsymbol{\cdot} xy^{\alpha_2} \boldsymbol{\cdot} xy^{\alpha_3} \mid S_3$. By changing the order in which the elements are obtained, we can show that all terms of $S_3$ are equal. In fact, 
$$xy^{\alpha_1} \cdot xy^{\alpha_2} = y^2 = xy^{\alpha_1} \cdot xy^{\alpha_3}$$ 
implies that $\alpha_2 \equiv \alpha_3 \pmod n$, and similarly $\alpha_1 \equiv \alpha_2 \equiv \alpha_3 \pmod n$. In particular, all terms of $S_3$ are equal. A similar result can be obtained if $|S_4| \ge 3$.

Since $n/2 \le k = |S_3|+|S_4|$, it follows that $\max\{|S_3|,|S_4|\} \ge n/4 \ge 2$. If $|S_3| = |S_4| = 2$ and $n=8$, then $L = 4 = n/2$, a contradiction. Therefore we assume that $\max\{|S_3|,|S_4|\} \ge 3$. If this maximum is $|S_4|$ and $xy^{\alpha} \mid S_4$ ($\alpha$ is even), then $xy^{\alpha} \cdot xy^{\alpha} = y^{2\alpha} \in C_{n/4} = \langle y^4 \rangle \not\ni y^2$, a contradiction. Therefore, $|S_4| \le 2$ and $|S_3| \ge 3$.

We have so far $|S_1| = n-k$, $|S_2| = 0$, $|S_3| \ge 3$, $|S_4| \le 2$ and $|S_3| + |S_4| = k$. Consider the values of $|S_4|$:

\begin{enumerate}[{(d.5.}i{)}]
\item 
Sub-subcase $|S_4| = 2$. Let $S_4 = xy^{\alpha_1} \boldsymbol{\cdot} xy^{\alpha_2}$, and let $(xy^{\beta})^{[2]} \mid S_3$ (that is, $\alpha_1$ and $\alpha_2$ are even, and $\beta$ is odd). It follows that 
$$y^{2\beta + n/2} = xy^{\beta} \cdot xy^{\beta} = y^2 = xy^{\alpha_1} \cdot xy^{\alpha_2} = y^{\alpha_1 + \alpha_2}.$$ 
We reorder only these four terms, obtaining the product 
$$\quad \quad \quad \quad xy^{\beta} \cdot xy^{\alpha_1} \cdot xy^{\beta} \cdot xy^{\alpha_2} = y^{\alpha_1+\beta+n/2} \cdot y^{\alpha_2+\beta+n/2} = y^{\alpha_1 + \alpha_2 + 2\beta} = y^{n/2 + 4},$$ 
and we complete the previous product with the subsequence $(y^2)^{[n/4-2]}$ over $C_{n/2}$, obtaining a product-one subsequence.

\item 
Sub-subcase $|S_4| = 1$. We may assume without loss of generality that $S = y^{[n-k]} \boldsymbol{\cdot} (xy^{\alpha})^{[k-1]} \boldsymbol{\cdot} xy^{\beta}$, where $\alpha$ is odd and $\beta$ is even. Since $y^2 = (xy^{\alpha})^2 = y^{2\alpha+n/2}$, we have that $\alpha \in \{n/4+1, 3n/4+1\}$. Suppose $\alpha = n/4+1$ (the case $\alpha = 3n/4+1$ is completely similar). If $k < n$, then it is possible to obtain both products 
$$\quad \quad \quad xy^{\beta} \cdot xy^{n/4+1} \cdot y = y^{n/4+\beta+2} \quad \text{ and } \quad xy^{n/4+1} \cdot xy^{\beta} \cdot y = y^{3n/4+\beta+2}.$$
Both exponents are even numbers and, modulo $n$, some of them is at least $n/2$, therefore it is possible to complete such product with some $y^2$'s, obtaining a product-one subsequence. It implies that $k = n$ and $S = (xy^{\alpha})^{[n-1]} \boldsymbol{\cdot} xy^{\beta}$, where $\alpha$ is odd and $\beta$ is even.

\item 
Sub-subcase $|S_4| = 0$. We may assume without loss of generality that $S = y^{[n-k]} \boldsymbol{\cdot} (xy^{n/4+1})^{[k]}$. If $k \le n-2$, then $$xy^{n/4+1} \cdot y \cdot xy^{n/4+1} \cdot y = y^{n/2 + 2},$$
therefore it is possible to complete such product with $(y^2)^{[n/4 - 1]}$, thus $S$ is not product-one free. If $k = n$, then $S = (xy^{n/4+1})^{[n]}$ and $(xy^{n/4+1})^n = 1$, thus $S$ is not product-one free. Hence, $k = n-1$.
\end{enumerate}
\end{enumerate}
\end{enumerate}
It completes the proof.

\qed

\end{document}